%% file: paper-arXiv.tex
\documentclass[12pt]{article}

\usepackage{amssymb,amsthm,amsmath}
\usepackage[mathscr]{eucal}

\input{declare}

\begin{document}
\author{Jay H. Beder \\ Margaret Ann McComack \\
Department of Mathematical Sciences\\University of Wisconsin-Milwaukee\\
P.O. Box 413\\Milwaukee, WI 53201-0413\\
\texttt{beder@uwm.edu}\\ \texttt{mamccomack@gmail.com} 
}

\title{A note on the minimum size of an orthogonal array}
\date{}

\maketitle

\vspace{-1cm}

\begin{abstract}
\input{abstract}
\end{abstract}


{\footnotesize {\bf AMS(MOS) subject classification.}
Primary: 62K15; 
Secondary:
05B15, 
62K05 
}

\section{Introduction}\label{intro-sec}
\input{intro}

\section{The numbers $L_t$} \label{Lt-sec}
\input{Lt}

\section{Uniqueness}  \label{uniqueness-sec}

\input{uniqueness}

\section{Constructions}  \label{constructions-sec}
\input{constructions}

\bibliographystyle{plain}

\end{document}

%% file: declare.tex
\usepackage{amssymb,amsthm,amsmath}
\usepackage[mathscr]{eucal}

\newcommand{\Z}{\mbox{$\mathbb Z$}}

\usepackage{bm}

\newcommand{\inv}{\ensuremath{^{-1}}}

\newcommand\set[1]{\{1, \ldots, #1 \}}
\DeclareMathOperator{\lcm}{lcm}
\DeclareMathOperator{\ord}{ord}
\newcommand{\co}{\ensuremath{\mspace{-5mu}:}}  
\newcommand{\cart}[3]{\ensuremath{{#1}_{#2} \times \cdots \times {#1}_{#3}}}

\newtheorem{theorem}{Theorem}[section]

\theoremstyle{definition}
\newtheorem{example}{Example}

%% file: abstract.tex
It is an elementary fact that the size of an orthogonal array of strength $t$ on $k$ factors must be a multiple of a certain number, say $L_t$, that depends on the orders of the factors.  Thus $L_t$ is a lower bound on the size of arrays of strength $t$ on those factors, and is no larger than $L_k$, the size of the complete factorial design.  We investigate the relationship between the numbers $L_t$, and two questions in particular:  For what $t$ is $L_t < L_k$?  And when $L_t = L_k$, is the complete factorial design the only array of that size and strength $t$?  Arrays are assumed to be mixed-level.

We refer to an array of size less than $L_k$ as a \emph{proper fraction}.  Guided by our main result, we construct a variety of mixed-level proper fractions of strength $k-1$ that also satisfy a certain group-theoretic condition. 

%% file: intro.tex
Let $D$ be an orthogonal array of size $N$ on $k$ factors, the $i$th factor having $s_i$ values or levels.  It is easy to see that if $t$ is the strength of $D$, then $N$ is a multiple of
  \begin{equation} \label{Lt-def-eq} L_t = \lcm\left( \prod_{i \in I} s_i,  |I| = t \right), \end{equation}
where $|I|$ is the number of indices in $I$.  These points are reviewed in detail below.

Obviously $L_t$ is a lower bound for $N$, although an array of strength $t$ and size $N = L_t$ may not exist.  Moreover, $L_t \le s_1 \cdots s_k = L_k$, the size of the complete factorial design on these factors, since the complete design certainly has strength~$t$.  From the point of view of applications, it is of interest to know when $L_t < L_k$, since in this case one may seek a design that is a proper fraction of the full factorial design.  Theorem~\ref{Lt-ineq-thm} gives a simple criterion based on the orders $s_i$ to determine those $t$ for which we have $L_t < L_k$.

\subsection{Definitions and notation}

As indicated above, we use $|I|$ to denote the cardinality of the set $I$.

A \emph{complete factorial design} is a Cartesian product $\cart{A}{1}{k}$, where  $A_i$ is a finite set.  In statistical design, $A_i$ is the set of \emph{levels} of factor $i$, and the elements of $\cart{A}{1}{k}$ are \emph{treatment combinations}.  We will refer to $s_i = |A_i|$ as the \emph{order} of factor $i$.  An \emph{orthogonal array} or \emph{design} $D$ on these factors is a multisubset
(a subset with possible repetitions) of $\cart{A}{1}{k}$. The \emph{size} of the array, $N$, is the number of elements (or \emph{runs}), counting multiplicities.  The design is \emph{symmetric} if $s_1 = \cdots = s_k$, and otherwise is asymmetric or \emph{mixed-level}.  We will say that a design is a \emph{proper fraction} if it is a  a proper subset of a complete factorial design.

If we write the elements of a design $D$ as columns, then we may represent $D$ as a $k \times N$ matrix of symbols.
The \emph{projection} of $D$ on $j$ factors is then the $j \times N$ submatrix consisting of the rows corresponding to those factors.  We say that $D$ has \emph{strength} $t$ if for any subset $I \subset \set{k}$ of size $t$, the projection of $D$ on the factors indexed by $I$ consists of $\lambda_I$ copies of the Cartesian product $\prod_{i \in I} A_i$, for some $\lambda_I$.  Evidently we have
  \[ N = \lambda_I \prod_{i \in I} s_i, \]
from which it follows easily that $N$ is a multiple of the number $L_t$ given by (\ref{Lt-def-eq}).   In particular, $L_1 = \lcm(s_1, \ldots, s_k)$ and $L_k = s_1 \cdots s_k$, and it is easy to see that
  \[ L_t | L_{t+1} \]
for each $t = 1, \ldots, k-1$ (each product in $L_t$ is contained in a product in $L_{t+1}$).  Therefore  \begin{equation} \label{weak-ineq} L_1 \le L_2 \le \cdots \le L_k. \end{equation}
In a symmetric design with $s_1 = \cdots = s_k = s$ we have $L_t = s^t$.

For each $t$, $L_t$ is the smallest possible size of an array of strength $t$.  An array of this size may not actually exist.  For example, in a symmetric design on $k=4$ factors, each with 6 levels, we have $L_2 = 36$, but there is no array of strength 2 and size 36, since this would be equivalent to two mutually orthogonal Latin squares of order 6 (see, e.g., \cite{Raghavarao}, pages 11 and 33).

As indicated above, it is of interest to know for which $t$ we have $L_t < L_k$.  Theorem~\ref{Lt-ineq-thm} gives a useful criterion, and also strengths the inequality (\ref{weak-ineq}).

%% file: Lt.tex
For each $I \subset \set{k}$ let  \begin{align*}
  e_I &= \gcd(s_i, i \in I),
\intertext{and consider those $I$ for which $e_I > 1$. Let}
  d &= \max\{ |I|: e_I > 1 \}. \end{align*}
Thus $|I| > d $ implies $e_I = 1$.  We see that $1 \le d \le k$, and $d = 1$ iff the orders $s_i$ are pairwise relatively prime. At the other extreme, $d=k$ iff $\gcd(s_1, \ldots, s_k) > 1$.  For example, $d=k$ if the design is symmetric.

In the proof below we need to choose a subset $I$ of size $d$ for which $e_I > 1$.  There can be more than one such subset: for example, if the numbers $s_i$ are $8, 12, 18, 27$, then $d = 3$, since the four numbers have no common factor but $\{8, 12, 18\}$ and $\{12, 18, 27\}$ have gcd greater than 1;  the corresponding sets $I \subset \set{4}$ are $\{1, 2, 3\}$ and $\{2, 3, 4\}$.  

\begin{theorem} We have $L_1 < \cdots < L_d = \cdots = L_k$.  In particular, $L_t = L_k$ iff $t \ge d$.
\label{Lt-ineq-thm} \end{theorem}

Before proving this theorem, we mention two useful facts. First, if $S_1, \ldots, S_n$ are sets of positive integers and $S = \cup_i S_i$, then
  \begin{equation} \label{lcm(lcm)-eq} \lcm(S) = \lcm(\lcm(S_1), \ldots, \lcm(S_n)). \end{equation}
Second, for positive integers $a_1, \ldots, a_n$ we have
  \begin{equation} \label{prod/gcd-eq} \lcm\left(\prod_{i \in I}a_1, \; |I| = n-1 \right) = \frac{a_1 \cdots a_n}{\gcd(a_1, \ldots, a_n)}. \end{equation}
These can both be proved prime-by-prime: Define $\ord_p(b)$, the \emph{order} of the prime $p$ in $b$, to be the largest power of $p$ dividing $b$.  To prove (\ref{prod/gcd-eq}), for example, let $f_i = \ord_p(a_i)$.  Then (\ref{prod/gcd-eq}) is the statement that for each prime $p$,
  \[ \max_j\left(\sum_{i=1}^n f_i - f_j\right) = \sum_{i=1}^n f_i - \min_j f_j. \]
(Property (\ref{prod/gcd-eq}) is familiar in the case $n=2$, namely, $ab = \lcm(a,b)\gcd(a,b)$.)

\begin{proof}[Proof of Theorem \ref{Lt-ineq-thm}]  Fix $t$. In order to compare $L_t$ and $L_{t+1}$, we begin by organizing the subsets $I \subset \set{k}$ of size $t$ into overlapping families $B_J$ indexed by the subsets $J$ of size $t+1$.  Namely, we put \begin{align}
B_J &= \{ I \subset J\co |I| = t \}.  \nonumber
\intertext{Let $B = \{ I \subset \set{k}\co |I| = t \}$.  Then}
   B &= \bigcup_{|J| = t+1} B_J. \nonumber
\intertext{Now from (\ref{lcm(lcm)-eq}) and (\ref{prod/gcd-eq}) we have}
   L_t &= \lcm\left( \prod_{i \in I} s_i,  |I| = t \right)  \nonumber \\
       &= \lcm \left( \lcm\left( \prod_{i \in I} s_i,  I \in B_J \right), |J| = t+1 \right) \nonumber \\
       &= \lcm \left( \prod_{i \in J} s_i/e_J, |J| = t+1 \right). \label{almost-Lt+1-eq}
   \end{align}

If $t \ge d$ then the expression (\ref{almost-Lt+1-eq}) is exactly $L_{t+1}$.  For if $t \ge d$ then $e_J = 1$ whenever $|J| = t+1$.  Thus $t \ge d$ implies that $L_t = L_{t+1}$.

If $t < d$, we claim that (\ref{almost-Lt+1-eq}) is strictly less than $L_{t+1}$.
To see this, fix a prime $p$ that divides exactly $d$ of the numbers $s_1, \ldots, s_k$.  (This is possible as there exists a set $I \subset \set{k}$ with $|I| = d$ and $e_{I} > 1$; let $p | e_{I}$.)
Our claim follows as soon as we show that the order of $p$ in (\ref{almost-Lt+1-eq}) is strictly less than $\ord_p(L_{t+1})$.  To do this, we need to show that this holds for each term in (\ref{almost-Lt+1-eq}), that is,
  \begin{equation} \label{order-ineq} \ord_p\left( \prod_{i \in J} s_i/e_J \right) < \ord_p(L_{t+1}) \end{equation}
for each $J\subset \set{k}$ with $|J| = t+1$.

Let $f_i = \ord_p(s_i)$.  Note that exactly $d$ of the integers $f_i$ are positive.  Renumbering if necessary, we may assume without loss of generality that
  \[ f_1 \ge f_2 \ge \cdots \ge f_d > 0 = f_{d+1} = \cdots = f_k. \]

Now on the one hand, the order of $p$ in $L_{t+1}$ is the maximum of the orders of $p$ in the products $\prod_{i \in J} s_i$, $|J| = t+1$.  But this maximum is attained for $J = \set{t+1}$, and equals $f_1 + \cdots + f_{t+1}$;  that is,
  \begin{equation} \label{ordp(Lt+1)-eq} \ord_p(L_{t+1}) =  \sum_{i=1}^{t+1} f_i.
\end{equation}

On the other hand, for any $J \subset \set{k}$ we have
  \[ \ord_p (e_J) = \min_{i \in J} f_i,   \]
and so
  \[ \ord_p \left( \prod_{i \in J} s_i / e_J \right) = \sum_{i \in J} f_i - \min_{i \in J} f_i. \]
We focus on sets $J$ with $|J| = t+1$, and consider two cases:  \begin{itemize}
  \item $\min_{i \in J} f_i > 0$: Here \[  \sum_{i \in J} f_i - \min_{i \in J} f_i <  \sum_{i \in J} f_i \le \sum_{i=1}^{t+1} f_i. \]
  \item $\min_{i \in J} f_i = 0$: In this case, $f_i = 0$ for at least one $i \in J$, so \[ \sum_{i \in J} f_i - \min_{i \in J} f_i = \sum_{i \in J} f_i  < \sum_{i = 1}^{t+1} f_i. \]
\end{itemize}
In either case (\ref{order-ineq}) holds, and therefore the order of $p$ in (\ref{almost-Lt+1-eq}) is less than (\ref{ordp(Lt+1)-eq}), which is what we needed to prove.
\end{proof}

Theorem \ref{Lt-ineq-thm} provides useful information about possible orthogonal arrays without our directly having to calculate the numbers $L_t$. In particular, it gives a necessary condition for the existence of a proper fraction of given strength, as illustrated by the examples below.
A special case of the theorem is given in \cite{Maggie-thesis}.
\medskip

\begin{example} In a $2^i 3^j$ experiment, that is, one having $i+j$ factors of which $i$ have 2 levels and $j$ have 3, we see that $d = \max(i,j)$.  In this case, no proper fraction has strength $t \ge \max(i,j)$.
\end{example}

\begin{example} Consider a $2 \times 3 \times 5 \times 6 \times 10 \times 15$ experiment.  Any subset of $ \{ 2, 3, 5, 6, 10, 15\}$ of size 4 must have $e_I = 1$, as two of its elements must be relatively prime.  On the other hand, there are subsets of size 3 that have $e_I > 1$ -- for example, $\{2, 6,10\}$ (that is, $I = \{1, 4, 5\}$).  Thus $d=3$, and so no proper fraction has strength 3.
\end{example}

\begin{example}  In order for there to exist a proper fraction of strength $k-1$ in an experiment with $k$ factors, it is necessary that $d=k$.  That is, the orders $s_1, \ldots, s_k$ must share a common factor.  In Section~\ref{constructions-sec} we give methods of construction that generate a number of proper fractions of strength $k-1$ that also satisfy a group-theoretic condition when the set of treatment combinations is a nonabelian group.  The condition and its origin are discussed there. \label{strength(k-1)-ex}
\end{example}

%
%

%% file: uniqueness.tex
When $L_t = L_k$, the smallest array of strength $t$ has the same size as the complete factorial design, and it is natural to ask whether the complete factorial is the only array of its size of strength $t$.  It is easy to see that for strength $t=1$ this is not true.  For suppose factor $i$ is indexed by the set $A_i$ where $|A_i| = s_i$.  Fill a $k \times L_k$ matrix by putting $L_k/s_i$ copies of $A_i$ in row $i$, in an arbitrary order.  The resulting array has strength $t=1$ and size $L_k$, but it is easy to fill the rows in such a way that the columns do not consist of all the elements of $\cart{A}{1}{k}$.

Of course, the rows of an array of strength 1 can be filled independently of each other, and so one might conjecture that such a construction is impossible if we require strength 2 or higher, but in fact this is not true either.  Consider a $3 \times 2 \times 2$ factorial experiment and let $A_1 = \{0, 1, 2\}$ and $A_2 = A_3 = \{0, 1\}$.  In this case $L_2 = 12$, the size of the complete factorial design, and so there is no smaller orthogonal array of size 12 and strength 2 on these factors.  We easily check that the $3 \times 12$ array
  \[ \left(
  \begin{array}{cccccccccccc}
    0 & 0 & 1 & 1 & 2 & 2 & 0 & 0 & 1 & 1 & 2 & 2 \\
    0 & 0 & 0 & 0 & 0 & 0 & 1 & 1 & 1 & 1 & 1 & 1 \\
    0 & 0 & 0 & 1 & 1 & 1 & 1 & 1 & 0 & 1 & 0 & 0 \\
  \end{array}
\right) \]
has strength 2, but it is clearly not the complete design, as some ordered triples occur more than once while others are missing.  It would be interesting to know whether there are applications in which such an array might be a useful alternative to the full factorial design. 

%% file: constructions.tex
In this section we construct a number of mixed-level orthogonal arrays of strength $k-1$ on $k$ factors.  They are all proper fractions and in addition satisfy a group-theoretic property (a ``conjugacy" condition) given in \cite{Beder&Beder14} that we describe below.  As we saw in Example~\ref{strength(k-1)-ex}, Theorem~\ref{Lt-ineq-thm} requires that the orders $s_i$ share a common factor.

Our examples consist of $k=3$ to 6 factors, with orders $s_i \le 10$.  Verification that each array satisfies our two conditions (strength $k-1$ and conjugacy) was done by computer, for which the code can be found in \cite{Maggie-thesis}.  Only those cases listed in Table~\ref{array-table} (below) have been checked.

\subsection{Conjugacy}
Throughout this section we will suppose that the sets $A_i$ indexing the levels of the factors are groups, and will use $G_i$ rather than $A_i$ as a more suggestive notation.  The set $G = \cart{G}{1}{k}$ of treatment combinations is a group, and is abelian iff all the factors $G_i$ are abelian.

Recall that the elements $a$ and $b$ of $G$ are conjugate if $b = gag\inv$ for some $g \in G$.  Conjugacy is an equivalence relation on $G$, and so partitions $G$ into \emph{conjugacy classes}.  In an abelian group, the conjugacy classes are all singleton sets.  It is an elementary fact that the conjugacy classes of $\cart{G}{1}{k}$ are of the form $\cart{C}{1}{k}$ where $C_i$ is a conjugacy class of $G_i$. The condition we require is that \emph{the design $D$ be a union of conjugacy classes of $G$}.\footnote{The condition assumed in \cite{Beder&Beder14} is that the counting function of the multiset $D$ be constant on conjugacy classes.  Since we are constructing proper fractions, $D$ is a set and its counting function is its ordinary indicator function, which is constant on conjugacy classes iff $D$ is a union of such classes.}

In all our examples $G_1$, and therefore $G$, is a nonabelian group, while $G_i$ will be abelian for $i >1$.  Thus the conjugacy classes of $G$ are essentially those of $G_1$ -- namely, of the form $C \times \{g_2\} \times \cdots \times \{g_k\}$ where $C$ is a conjugacy class of $G_1$ and $g_i \in G_i$. The nonabelian groups we will use are the following, where $e$ will denote the identity element.
\begin{itemize}
\item $S_3$, the symmetric group on 3 letters.  We have $S_3 = \{e, x, y, a, b, c\}$, where $x$ and $y$ are the 3-cycles and $a, b$ and $c$ are the 2-cycles.
    We will use either of two orderings:  \begin{enumerate}
    \item $e \; | \;x \;y \;|\; a \;b \;c$.
    \item $e \;|\; a \;b \;c\; | \;x \;y$.
    \end{enumerate}
    The vertical lines are inserted merely to indicate the conjugacy classes.
\item  $Dih_4 = \{e, q, r, s, a, b, x, y \}$, the symmetries of the square, where $q$ is the half-turn, $r$ and $s$ are the quarter-turns, $a$ and $b$ are the reflections about the diagonals, and $x$ and $y$ are the reflections about the lines joining the midpoints of the opposite sides.  The conjugacy classes are $ e\; |\; q\; |\; r\; s\; |\; a\; b\; |\; x\; y $.
\item $Dih_5 = \{e, a, b, c, d, v, w, x, y, z\}$, the symmetries of the regular pentagon, where $a$ and $d$ are the rotations of $\pm 72^{\circ}$, $b$ and $c$ are the rotations of $\pm 144^{\circ}$, and $v, w, x, y$ and $z$ are the reflections about a line from a vertex to the midpoint of the opposite side.  The conjugacy classes are $e\; |\; a\; d\; |\; b\; c\; |\; v\; w\; x\; y\; z$.
\end{itemize}
In each of the above, the ordering of conjugacy classes is fixed, and the ordering of elements within each class is fixed but arbitrary.  These fixed orderings are assumed in the constructions below. In general, the group $Dih_n$ is a \emph{dihedral group}, and represents the symmetries of the regular $n$-gon. Of course, $S_3$ is also $Dih_3$.

For abelian groups we will use $\Z_n$, the integers modulo $n$, whose elements are $0, 1, \ldots, n$.  We will fix this order.

\subsection{Construction}

If $s_1 = 6, 8$ or 10, we let $G_1 = S_3, Dih_4$ or $Dih_5$, respectively.  For $i \ge 2$, we let $G_i = \Z_{s_i}$.  Our construction has the following cases.

\textbf{gcd($s_1, \ldots, s_k$) = 2 or 4}:  We construct an array of minimal size $L_{k-1}$.  Defining the integers  \begin{align*}
v_1 &= \frac{L_{k-1}}{s_1},  \\
v_j &= \frac{L_{k-1}}{s_2\cdots s_j} \; \mbox{for} \; j \ge 2,
\end{align*}
we fill in the rows of the array as follows:\medskip

Row 1: Write the elements of $G_1$ in the fixed order, the whole sequence repeated $v_1$ times.  If $s_1 = 6$, use $S_3$ in its first ordering.

Row 2: Write the elements of $G_2 (= \Z_{s_2})$ in the fixed order, repeating each element $v_2$ times.

Rows 3 through $k-1$: Write the elements of $G_j (= \Z_{s_j})$ in the fixed order, repeating each element $v_j$ times.  Repeat the whole pattern until the row is complete (a total of $L_{k-1}$ entries).

Row $k$: Write the elements of $G_k$, repeating each element $v_k$ times, in the following pattern:  first in the given order, then in the reverse order, alternating in that way until the row is complete.
\medskip

\begin{example}  The following arrays illustrate the construction method we have described.  \begin{itemize}

\setlength{\arraycolsep}{1mm}

\item 1/2 fraction of a $6 \times 2 \times 2 \times 2$ design, strength 3.  Here $G = S_3 \times \Z_2 \times \Z_2 \times \Z_2$.
\[\left(
        \begin{array}{cccccccccccccccccccccccc}
          e & x & y & a & b & c & e & x & y & a & b & c & e & x & y & a & b & c & e & x & y & a & b & c \\
          0 & 0 & 0 & 0 & 0 & 0 & 0 & 0 & 0 & 0 & 0 & 0 & 1 & 1 & 1 & 1 & 1 & 1 & 1 & 1 & 1 & 1 & 1 & 1 \\
          0 & 0 & 0 & 0 & 0 & 0 & 1 & 1 & 1 & 1 & 1 & 1 & 0 & 0 & 0 & 0 & 0 & 0 & 1 & 1 & 1 & 1 & 1 & 1 \\
          0 & 0 & 0 & 1 & 1 & 1 & 1 & 1 & 1 & 0 & 0 & 0 & 0 & 0 & 0 & 1 & 1 & 1 & 1 & 1 & 1 & 0 & 0 & 0 \\
        \end{array}
\right)\]

\item A 1/2 fraction of an $8 \times 2 \times 2$ design, strength 2.
$G = Dih_4 \times \Z_2 \times \Z_2$.

\[\left(
        \begin{array}{cccccccccccccccc}
          e & q & r & s & a & b & x & y & e & q & r & s & a & b & x & y \\
          0 & 0 & 0 & 0 & 0 & 0 & 0 & 0 & 1 & 1 & 1 & 1 & 1 & 1 & 1 & 1 \\
          0 & 0 & 0 & 0 & 1 & 1 & 1 & 1 & 1 & 1 & 1 & 1 & 0 & 0 & 0 & 0 \\
        \end{array}
\right)\]

\item A 1/2 fraction of a $10 \times 2 \times 2$ design, strength 2.  $G = Dih_5 \times \Z_2 \times \Z_2$

\[\left(
        \begin{array}{cccccccccccccccccccc}
          e & a & b & c & d & v & w & x & y & z & e & a & b & c & d & v & w & x & y & z \\
          0 & 0 & 0 & 0 & 0 & 0 & 0 & 0 & 0 & 0 & 1 & 1 & 1 & 1 & 1 & 1 & 1 & 1 & 1 & 1 \\
          0 & 0 & 0 & 0 & 0 & 1 & 1 & 1 & 1 & 1 & 1 & 1 & 1 & 1 & 1 & 0 & 0 & 0 & 0 & 0 \\
        \end{array}
\right)\]

\item A 1/4 fraction of an $8 \times 4 \times 4$ design, strength 2.  $G = Dih_4 \times \Z_4 \times \Z_4$.  Vertical lines are inserted to reveal the pattern in the last row.

\[\left(
        \begin{array}{cccccccc|cccccccc|cccccccc|cccccccc}
        e & q & r & s & a & b & x & y & e & q & r & s & a & b & x & y & e & q & r & s & a & b & x & y & e & q & r & s & a & b & x & y\\
        0 & 0 & 0 & 0 & 0 & 0 & 0 & 0 & 1 & 1 & 1 & 1 & 1 & 1 & 1 & 1 & 2 & 2 & 2 & 2 & 2 & 2 & 2 & 2 & 3 & 3 & 3 & 3 & 3 & 3 & 3 & 3\\
        0 & 0 & 1 & 1 & 2 & 2 & 3 & 3 & 3 & 3 & 2 & 2 & 1 & 1 & 0 & 0 & 0 & 0 & 1 & 1 & 2 & 2 & 3 & 3 & 3 & 3 & 2 & 2 & 1 & 1 & 0 & 0\\
        \end{array}
\right)\]
\end{itemize}

\end{example}
\medskip

\textbf{gcd($s_1, \ldots, s_k$) = 6}:  We in fact consider only symmetric ``$6^k$" designs, so that $L_{k-1} = 6^{k-1}$.  The method described above may be used to create arrays of this minimum size (and therefore 1/6 fractions), but the arrays will not satisfy the conjugacy requirement.  If we modify the integers $v_i$ as follows: \begin{align*}
v_1 &= \frac{3L_{k-1}}{s_1} = 3 \cdot 6^{k-2},  \\
v_j &= \frac{3L_{k-1}}{s_2\cdots s_j} = 3\cdot 6^{k-j} \; \mbox{for} \; j \ge 2,
\end{align*}
then the method will produce 1/2 fractions of strength $k-1$.  Note that the first $v_k = 3$ elements of $G_1 = S_3$ is a union of conjugacy classes of $S_3$.

We don't present an example since $N = (1/2)6^k$ rather large.
\medskip

\textbf{gcd($s_1, \ldots, s_k$) = 3}:  We apply this to $6 \times 3 \times 3 \cdots$ factorial experiments.  Here we alter both the integers $v_i$ and the steps of construction, since the original steps will produce arrays of minimum size (1/3 fractions) that do not satisfy the conjugacy property.  We set \begin{align*}
v_1 &= \frac{2L_{k-1}}{s_1} = 2\cdot 3^{k-2},  \\
v_j &= \frac{2L_{k-1}}{s_2\cdots s_j}= 4\cdot 3^{k-j} \; \mbox{for} \; j \ge 2,
\end{align*}
and fill in the rows of the array as follows:\medskip

Row 1: We use the \emph{second} order of $S_3$. Write the elements of $G_1 = S_3$ in this fixed order, then in the reverse order, alternating until there are $N=2L_{k-1}$.

Rows 2 through $k-1$: These steps are identical to those given above.

Row $k$: Write the elements of $G_k = \Z_3$, repeating each element $v_k = 4$ times.  Then do the same, but permuting the elements of $\Z_3$ cyclically, and again with another cyclic permutation, continuing in this pattern until the row is filled.

This method produces 2/3 fractions.
\medskip

\begin{example} A 2/3 fraction of a $6 \times 3 \times 3$ design, strength 2. Here $G = S_3 \times \Z_3 \times \Z_3$. Vertical lines are inserted to reveal the pattern in the last row.

\setlength{\arraycolsep}{1mm}
\[\left(
        \begin{array}{cccccccccccc|cccccccccccc|cccccccccccc}
          e & a & b & c & x & y & y & x & c & b & a & e & e & a & b & c & x & y & y & x & c & b & a & e & e & a & b & c & x & y & y & x & c & b & a & e\\
          0 & 0 & 0 & 0 & 0 & 0 & 0 & 0 & 0 & 0 & 0 & 0 & 1 & 1 & 1 & 1 & 1 & 1 & 1 & 1 & 1 & 1 & 1 & 1 & 2 & 2 & 2 & 2 & 2 & 2 & 2 & 2 & 2 & 2 & 2 & 2 \\
          0 & 0 & 0 & 0 & 1 & 1 & 1 & 1 & 2 & 2 & 2 & 2 & 2 & 2 & 2 & 2 & 0 & 0 & 0 & 0 & 1 & 1 & 1 & 1 & 1 & 1 & 1 & 1 & 2 & 2 & 2 & 2 & 0 & 0 & 0 & 0 \\
        \end{array}
\right)\]
\end{example}
\medskip

Table \ref{array-table} summarizes 31 arrays constructed using the methods we have described.

\begin{table}
\caption{This table lists arrays constructed using the methods in this paper.  Strength is $k-1$ where $k$ is the number of factors.  All arrays are if minimal size $L_{k-1}$ unless otherwise indicated. Examples are separated into four groups according to $\gcd(s_1, \ldots, s_k)$.} \label{array-table}
\centering
\begin{tabular}{c c c c c}
\hline\hline
Complete & Size of         & Size of & &   \\
Design   & Complete Design & Array &   & Fraction \\[0.5ex]
\hline
$6 \times 2 \times 2^*$ & 24 & 12 & & $1/2$ \\
$6 \times 2 \times 2 \times 2^*$ & 48 &  24 & & $1/2$ \\
$6 \times 4 \times 4$ & 96 &  48 & & $1/2$ \\
$6 \times 4 \times 4 \times 4$ & 384 &  192 & & $1/2$ \\
$6 \times 4 \times 2$ & 48 &  24 & & $1/2$ \\
$6 \times 6 \times 2$ & 72 &  36 & & $1/2$ \\
$6 \times 6 \times 4$ & 144 &  72 & & $1/2$ \\
$8 \times 2 \times 2$ & 32 &  16 & & $1/2$ \\
$8 \times 2 \times 2 \times 2$ & 64 &  32 & & $1/2$ \\
$8 \times 2 \times 2 \times 2 \times 2$ & 128 &  64 & & $1/2$ \\
$8 \times 2 \times 2 \times 2 \times 2 \times 2$ & 256 &  128 & & $1/2$ \\
$8 \times 6 \times 6$ & 288 & 144 & & $1/2$ \\
$8 \times 6 \times 6 \times 6$ & 1,728 &  864 & & $1/2$ \\
$8 \times 4 \times 2$ & 64 &  32 & & $1/2$ \\
$8 \times 6 \times 2$ & 96 &  48 & & $1/2$ \\
$8 \times 6 \times 4$ & 192 &  96 & & $1/2$ \\
$10 \times 2 \times 2$ & 40 &  20 & & $1/2$ \\
$10 \times 2 \times 2 \times 2$ & 80 &  40 & & $1/2$ \\
$10 \times 4 \times 4$ & 160 &  80 & & $1/2$ \\
$10 \times 4 \times 4 \times 4$ & 640 &  320 & & $1/2$ \\
$10 \times 6 \times 6$ & 360 &  180 & & $1/2$ \\
$10 \times 6 \times 6 \times 6$ & 2,160 &  1,080 & & $1/2$ \\
$10 \times 4 \times 2$ & 80 &  40 & & $1/2$ \\
$10 \times 6 \times 2$ & 120 &  60 & & $1/2$ \\
$10 \times 6 \times 4$ & 240 &  120 & & $1/2$ \\ \hline
$8 \times 4 \times 4$ & 128 &  32 & & $1/4$ \\
$8 \times 4 \times 4 \times 4$ & 512 &  128 & & $1/4$ \\ \hline
$6 \times 6 \times 6$ & 216 &  108 & $=3L_{k-1}$ & $1/2$ \\
$6 \times 6 \times 6 \times 6$ & 1,296 &  648 & $=3L_{k-1}$ & $1/2$ \\ \hline
$6 \times 3 \times 3$ & 54 &  36 & $=2L_{k-1}$ & $2/3$ \\
$6 \times 3 \times 3 \times 3$ & 162 &  108 & $=2L_{k-1}$ & $2/3$ \\
\hline
\multicolumn{4}{c}{* Another example is constructed in \cite{Beder&Beder14}.}
\end{tabular}
\end{table} 

%% file: paper-arXiv.bbl
\begin{thebibliography}{1}

\bibitem{Beder&Beder14}
Jay~H. Beder and Jesse~S. Beder.
\newblock Generalized wordlength patterns and strength.
\newblock {\em Journal of Statistical Planning and Inference}, 144:41--46,
  2014.

\bibitem{Maggie-thesis}
Margaret~Ann McComack.
\newblock Constructing orthogonal arrays on non-abelian groups.
\newblock Master's thesis, University of Wisconsin -- Milwaukee, 2013.

\bibitem{Raghavarao}
Damaraju Raghavarao.
\newblock {\em Constructions and Combinatorial Problems in Design of
  Experiments}.
\newblock Dover, 1988.
\newblock Reprint of the original Wiley edition, 1971.

\end{thebibliography}
